\documentclass{conm-p-l}

\newtheorem{theorem}{Theorem}[section]

\newtheorem{prop}[theorem]{proposition}
\theoremstyle{definition}

\theoremstyle{remark}

\numberwithin{equation}{section}

\def\tr{{\rm tr}}
\def\fs{+_{\mathcal{ F}}}
\def\gs{*_{\mathcal {G}}}

\begin{document}

\title{Informative Words and Discreteness}

\author{Jane Gilman }
\address{Department of Mathematics, Rutgers University, Newark, N.J. 07102}
\email{gilman@andromeda.rutgers.edu}
\thanks{This work has been supported in part by
NSA grant \#02G-186 and the Rutgers University Research
Foundation}

\date{}

\dedicatory{Dedicated to Gerhard Rosenberger on his $60^{th}$
birthday.}


\begin{abstract}

  There are certain families of words and word
sequences (words in the generators of a two-generator group) that
arise frequently in the Teichm{\"u}ller theory of hyperbolic
three-manifolds and  Kleinian and Fuchsian groups and  in the
discreteness problem for two generator matrix groups. We survey some
of the families of such words and sequences: the  semigroup of so
called {\sl good} words of Gehring-Martin, the so called {\sl
killer} words of Gabai-Meyerhoff-NThurston, the Farey words of
Keen-Series and Minsky,  the discreteness-algorithm Fibonacci
sequences of Gilman-Jiang and {\sl parabolic dust} words. We survey
connections between the families and establish a new connection
between good words and Farey words.
\end{abstract}

\maketitle

In 1970's and 80's Gerhard Rosenberger and his collaborators
developed a theory of discreteness for two generator subgroups of
$PSL(2,\mathbb{R})$. Essential to the theory was the concept of
replacing pairs of successive generators by Nielsen equivalent pairs
in a {\sl trace minimizing} manner (in order to reach a pair where a
discreteness determination could be made). Although they did not use
the word {\sl algorithm}, algorithms and sequences of words were
central to their conceptualization of the discreteness problem. In
this paper we survey a number of areas where their concepts have had
an impact upon the discreteness problem both directly or indirectly.
A partial list of the papers they wrote are included in the
bibliography.

\section{Introduction}

 Classical Riemann
surface theory is the basis of the development of a lot of modern
mathematics. In the last few decades greater interest has grown in
the theory of hyperbolic geometry and three-manifolds. Instead of
studying the surfaces or  hyperbolic thee-manifolds, one can
equivalently study the theory of their uniformizing groups, Fuchsian
groups and Kleinian groups, respectively. Because questions  about
these discrete matrix groups are often intractable, families of
surfaces or three-manifolds or their groups are studied by studying
their Teichm{\"u}ller spaces or their representation spaces or
parameter spaces. Recent important results about (appropriately
defined) parameter spaces and their boundaries have been obtained by
Gabai-Meyerhoff-N.Thurston, Gehring-Martin,
 Keen-Series-Maskit, and Minsky.
(See \cite{GMT},  \cite{GM3}, \cite{GM9},  \cite{KS1},
\cite{GMgood}, \cite{YM}, and the references given there.) During
this period, algebraic geometers and researchers  in symbolic
computation have been working on creating user friendly programs for
computing with on algebraic curves. (See the work of Buser and
Sepp{\"a}l{\"a} \cite{Bu}, \cite{BS1},  \cite{BSil}, and the
references given there.)
 A  further
thread has been the development of algorithms for determining the
discreteness or non-discreteness of two-generator subgroups of real
matrix groups. (See the work of Gilman, Gilman-Maskit and Jiang in
addition to the papers of Purzitsky and Rosenberger   \cite{G19},
\cite{G21}, \cite{GiMa},  \cite{YC}, \cite{PR}, \cite{PR2}, and
\cite{Ro}.)
 Here one of the main results is  the existence
of  algorithms that  solve the $PSL(2, \mathbb{R})$ discreteness
problem. The algorithms can
 be given in several different forms. In some cases they
can be implemented on a computer if,  for example, the entries in
the matrix are assumed to be algebraic numbers lying in a finite
extension of the rationals given in terms of their minimal
polynomials \cite{G19}. It can also be shown that (when
appropriately defined) the discreteness problem can not be solved
without an algorithm \cite{G21}.

 We view the results of Gerhing-Martin,
Gilman-Maskit-Jiang, Keen-Series, Gabai-Meyerhoff-N.Thurston,
Minsky, et al as being connected to the  results of Buser,
Sepp{\"a}l{\"a}, et al and their project of computing on algebraic
Riemann surfaces (CARS) and algorithm questions through the
observation that all are studying families of words in two-generator
groups and that these are precisely the families of words that any
algorithm in $PSL(2,\mathbb{C})$ would require.

The families of words studied include {\sl Farey  words}, so called
{\sl good words}, so called {\sl killer words}, Fibonacci and
non-Fibonacci word sequences, the algorithmic words and sequences of
pairs primitive associates in free groups. Precise definitions are
given in section \ref{section:wd}.

A current  goal is to put all of the different families of words
into a unified algebraic context. That is, is to establish
connections between the different families of words and to apply
these connections to obtain further results in the theory of
discrete groups, especially for $PSL(2, \mathbb{C})$ and results
applicable to computing on three-manifolds and algebraic curves. For
example, to address such problems as finding short geodesics on
Riemann surfaces or  finding a discreteness algorithms in the case
of complex matrix entries.

In this paper we survey some of the families of words that are
currently in use and known connections among them. We also establish
some new connections. The words begin with Rosenberger and the
concepts of replacing generators by Nielsen equivalent generators in
a {\sl trace minimizing manner}. The replacement in a trace
minimizing manner is an algorithm. Both the algorithmic nature of
this problem and the trace minimizing are significant.

\section{Overview}

Let $\mathbb{M}$ denote the group of all M{\"o}bius transformations
of the extended complex plane $\overline{\mathbb{C}}= {\mathbb{C}}
\cup \{\infty\}$. We associate  with a M{\"o}bius transformation
\[f= \frac{az+b}{cz+d} \in \mathbb{M}, \\\ ad-bc =  1,\]
one of the two matrices that induce the action of $f$ on the
extended plane
\[A =  \left( \begin{array}{cc} a  & b \\ c & d \\
\end{array} \right) \in SL(2,{\mathbb{C}})\]
and set ${\tr}(f)= {\tr}(A)$ where ${\tr}(A)$ denotes the trace of
$A$.  For a two generator group,  once a  matrix corresponding to
each  generator is  chosen,  the matrices  associated to all other
elements of the group are determined and thus their traces are
determined. Elements of $\mathbb{M}$ are classified as loxodromic,
elliptic, hyperbolic or parabolic according to the nature of the
square of their traces.

A subgroup $G$  of  $\mathbb{M}$ is  discrete if  no  sequence of
distinct elements of $G$  converges to the identity. A discrete
non-elementary subgroup of  $\mathbb{M}$ is a {\sl Kleinian} group.
If the matrix entries are real (or if the group is conjugate to such
a group), then such a discrete group is called a {\sl Fuchsian}
group.

The problem of determining discreteness  for even a two-generator
real matrix group is non-trivial. It turns out, that if $G = <A,B>$
is generated by $A$ and $B$ in $\mathbb{M}$, an algorithm is needed
to determine discreteness \cite{G21}.  The discreteness algorithm
inputs $A$ and $B$ and outputs one of two statements: the group $G$
generated by the matrices is {\sl discrete} or the group is {\sl not
discrete}. The algorithm appears in a number of different forms. A
corresponding algorithm for two matrices in $PSL(2, \mathbb{C})$ is
not known to exist.

\section{Informative Words \label{section:wd}}

 We summarize types  words that
 have been
used to obtain significant results about Kleinian groups and their
representation space(s).

\subsection{Primitive words}
 Let $W=W(A,B)$ denote a word in the generators $A$ and
$B$ of $G = <A,B>$ so that $$W(A,B) = A^{u_1}B^{v_1}A^{u_2} \dots
B^{v_t} \mbox{  for some integers } u_1,...,u_t;  v_1, ...,
v_t.\;\;\;(*)$$

We are primarily interested in {\sl primitive} words, words that can
be extended to  a minimal set of generators for the group. In the
case of   primitive words that generate a two-generator free group,
 we may assume that up to cyclic permutation and interchanging the
 generators there is a unique shortest word where the $u_i$ are all equal to each other
 and are all $\pm 1$.
That is, we can apply results from  \cite{MKS} to see that if $G=
<A, B>$ is a free group on two generators, then  a primitive word
can be written in a unique canonical  form
$A^{-1}B^{v_1}A^{-1}B^{v_2}A \cdots $ or the equivalent with $A$
replacing $A^{-1}$  and/or $A$ and $B$ interchanged. We call the
$v_i$ the {\sl primitive exponents}.

\subsection{Algorithmic words}
In the geometric algorithm \cite{G19,G21} if discreteness or
non-discreteness cannot be determined directly from the pair of
generators, then the pair  is replaced by a new pair. At step  $k$,
the generating pair $A_k,B_k$ is replaced by a new pair $A_{k+1},
B_{k+1}$. However, one of the new generators is the same as one of
the previous generator. Thus the sequence of the  {\sl algorithm
words} refers to the sequence of new generators.

The fact that the replacement procedure stops to give an algorithm
as opposed to a procedure (something that does not necessarily stop)
depends upon the fact that the algorithm is trace minimizing.

There are two types of steps that
 occur in the algorithm, steps that increase the word length in a
linear manner and steps that make the word length increase as a
Fibonacci sequence \cite{G21}.  They are termed
Fibonacci/non-Fibonacci steps. The size of the input to the
algorithm is measured by the maximal initial trace.  In the analysis
of the algorithm, Jiang shows that when the word length increases
exponentially in a Fibonacci step, the absolute value of trace of
the word (\cite{YC}) decreases logarithmically. For non-Fibonacci
steps both are linear. This is key to proving that the algorithm has
polynomial complexity.

\subsection{Farey Words}
 Farey words have been used by Keen and Series \cite{KS1} to find
 the boundary of the Riley slice of Schottky space. Farey words also play a
 role in recent work of Minsky \cite{YM}.

 {\sl  Farey words} (apres the mathematician Farey)
 are related to continued fraction expansions and
 the tessellation generated by the action of
$PSL(2,\mathbb{Z})$ on the upper half plane  \cite{Ser1}. If $p/q$
is a rational
 number between $0$ and $1$
and if  $[a_0,...,a_k]$ is the continued faction expansion so that
${{{p}\over{ q}}}  = {1 \over {a_0 + {1
 \over {a_1 +{1 \over a_2 + \dots }}}}}$,  then there is a word,
denoted by $W_{p/q}$ which is the {\sl Farey word} corresponding to
$[a_0,...,a_k]$. Further it is shown in \cite{KS1} that there is
another sequence
 of integers $[v_1,...,v_t]$ obtained from $[a_0,...,a_k]$
so that $W_{p/q}$ is equal to the word $W(A,B)$ with primitive
exponents $[v_1,.....,v_p]$ where $t=p$, $u_i = -1 \forall i$ and
$\Sigma v_i=q$. Farey words that are neighbors, (i.e. that satisfy
$|ps-rq|=1$), can be added with the following Farey addition (when
${\frac{p}{q}} < {\frac{r}{s}}$): $W_{\frac{p}{q}} \fs
W_{\frac{r}{s}} = W_{\frac{p+r}{q+s}}$.

 We  write
${\frac{p}{q}} =[a_0,a_1,...,a_k]$ to denote the fact that
${\frac{p}{q}}$ has continued fraction expansion
 $[a_0,a_1,...,a_k]$.

\subsection{Trace Polynomials}
We recall the notion of a {\sl trace polynomial}. The trace of a
word $W=W(A,B)$  can be described by a polynomial in $\tr(A)$,
$\tr(B)$ and $\tr([A,B])$ or $\tr(AB)$  where $[A,B]$ is the
multiplicative commutator of $A$ and $B$. For example,  for $h \in
\mathbb{M}$ as above,  $\tr (h^2)$ is given by the polynomial $p(x)
= x^2 -2$ since $\tr ( h^2) = (\tr h)^2 - 2$. So families of words
give rise to families of so called trace polynomials, polynomials in
the three variables: $\tr A$, $\tr B$, $\tr [A,B]$. In addition to
using different families of words, different families of parameters
are also used in different settings. Thus the trace polynomials can
also be viewed as polynomials in the variables $\beta(f) =
(\tr^2f-4)$ and $\gamma(f,g) = \tr [f,g] -2$.

 In the case that  $G$ is discrete and Fuchsian,
then the axis of a hyperbolic  element of $G$ projects onto a
geodesic on the quotient surface and $\tr(A)$ gives  the length $T$
of this geodesic
 (specifically,  $|\tr A | = 2 \cosh ({1 \over 2}T)$). Thus trace
polynomials are essentially lengths of geodesics.  A similar
statement holds for the complex length of a loxodromic element.

\subsection{Good words}
{\sl Good Words}  were defined by  Gehring and Martin and used  to
obtain results about the $\gamma \beta$ parameter space and minimal
volume three-manifolds. If $G = <A,B>$ is the free group on $A$ and
$B$, then a {\sl good} word
 is a word in $A$ and $B$ that starts and ends with a power of $A$ or
 its inverse and the exponents of the $A$'s oscillate in sign.
Thus using the notation in formula $*$, $u_1 = \pm 1$ and $u_{i+1} =
-u_i$ and $g_i=v_i$, $v_i \ne 0 \mbox{ if } i\ne t$,  for example, a
good word can be written as $A^{-1}B^{g_1}AB^{g_2}A^{-1}B^{g_3}
\cdots A^{-1}$ \cite{GMgood}.

We let ${\mathcal W}$ denote  the family of {\sl good} words. These
good words give rise to a family of trace polynomials, $\{
{{\mathcal P}_{Wgood}}  \}.$
 Good words
 can be composed and this  gives ${\mathcal W}$ a semi-group structure.
If $W_1(A,B)$ and $W_2(A,B)$ are good words, then $W_1 \gs W_2 =
 W_1(W_2(A,B), B)$. That is $W_2(A,B)$ is substituted into $W_1(A,B)$
wherever $A$ occurs in $W_1$. To be more precise

\begin{theorem} [Gehring-Martin] Given $W \in {\mathcal W}$, then  $ \exists  P_W$, a
polynomial with integer coefficients, such that
 if $\rho: F_2 \rightarrow \mathbb{M}$
with $f= \rho(a), g= \rho(b) , h = \rho(w)$  with $\gamma = tr[f,g]
-2, \beta= \beta(f)$, then $P_W(\gamma, \beta) = \gamma(f,h)$.
\end{theorem}
Thus  $P_{W_1 \gs W_2}(\gamma, \beta) = P_{W_1}(P_{W_2}(\gamma,
\beta), \beta)$. Since  for any $f$ and $g$ in $\mathbb{M}$ one can
regard the trace parameters, as being $\beta(f) = (\tr f)^2 -4$,
$\beta(g) = (\tr g)^2 -4$
 and $\gamma(f,g)
 = \tr [f,g] - 2$,
This gives a new $\gamma$ parameter, a new way of  moving around
parameter space and it  differs from the traditional way of
producing new $\gamma$'s.

\subsection{Killer Words}
{\sl Killer} words arise in  \cite{GMT}, historically the first
paper to prove a result about three-manifold theory using computer
implementations. That work builds  upon ideas initially developed by
Riley  \cite{Ri}. The technique in \cite{GMT} is to divide the
potential representation space into grids. A {\sl killer} word is a
word in the generators whose trace becomes very small on a given
grid.  J{\o}rgensen's inequality shows that a group is not discrete
if the trace(s) of some  word(s) are  small enough \cite{J1}. A
procedure is used to  find killer words for appropriate grids, thus
ruling out points in a grid possessing a killer word  as potential
parameters for a discrete group. An exact description of the {\sl
killer} words is not really known for  not all of the killer words
used in \cite{GMT} come from a fixed algorithm. The program produces
some by trial and error. It would be nice to bring a theory to the
killer words. Gaven Martin and T. Marshall are in the process of
finding sets of killer words that are good words.

\subsection{The question} The main question is
\begin{quote}
Can one establish connections between all of these disparate
families of informative words?
\end{quote}
Typical of the type of connection we hope to make is the following.
We observe that a good word in $A$ and $B$ is a Farey word $W(C,D)$
in the generators of the group $G_0$ where $C=A$, $D = BAB^{-1}$ and
$G_0 = <A, BAB^{-1}>$. Thus results applied to Farey words in $G$
and the boundary of its representation space can be applied to the
Farey words in $G_0$ and the boundary of the representation space of
$G_0$ and vice-verse. A related question is, \begin{quote}Is there a
corresponding connection   between $\fs$, addition of Farey
neighbors,   and $\gs$  good word multiplication? \end{quote}

\section{Connections}
As noted earlier G. Martin and T.Marshall are working on connections
between good words and killer words. Connections between algorithmic
words and Farey words and between short geodesics and  the
complexity of the algorithm were found in \cite{GK}. They are
summarized below after some additional notation and terminology are
introduced. Connections between good words and Farey words are
mentioned. More details of the latter will appear in \cite{Gconn}.

\subsection{The $F$-sequence}
The algorithm begins with a
 pair of generators  $(A,B)$ (which have been appropriately
 normalized  by interchanging $A$ and $B$ and replacing $A$ and $B$ by $A^{\pm 1}$ and $B^{\pm 1}$ as necessary.
If it does not stop and say "$G$ is discrete" or "$G$ is not
discrete", it  determines {\sl the} next pair of generators, which
is either $((A^{-1}B)^{-1}, B)$ or $(B^{-1}, A^{-1}B)$. The former
type of replacement is called a linear or non-Fibonacci step and the
latter a Fibonacci step (see \cite{G21} and \cite{YC}). The
algorithm dictates which type of step is the next step. Thus one can
associate to the algorithm its {\sl $F$-sequence} the numbers
$[n_1,n_2,....,n_k]$ where $n_i$ denotes the number of consecutive
linear steps before the next
 Fibonacci step  (see \cite{GK}).

If we have a primitive word in $A$ and $B$ given by its F-sequence,
we can multiply the word out and write it in terms of its primitive
exponents, the  $v_i$. For example, we can expand the word
$(A^{-1}B^{n_1})\cdot[B\cdot ((A^{-1}\cdot B^{n_1})^{n_2})]^{n_3})
\cdots$  in the form $B^{v_0}A^{-1}B^{v_1}A^{-1}B^{v_2}A^{-1}B^{v_3}
\cdots B^{v_w}$ for some integer $w$. The the primitive exponents
can be determined by the F-sequence and vice-versa (\cite{GK}).

\subsection{Continued fraction expansions and the  left-right sequence}
A Farey word, $W_{p/q}=W(A,B)$, corresponds to a geodesic connecting
the point at infinity to the rational point $p/q$ on the boundary of
the upper-half-plane. It will cross successive triangles in the
 $PSL(2,\mathbb{Z})$ tessellation. After replacing the tessellation by that of
 an appropriate subgroup of $PSL(2,\mathbb{Z})$, such a  geodesic when exiting a
 given triangle, either
 enters the next left triangle or the next right triangle and thus corresponds to a
  sequence of  lefts and rights, denoted by
 $L^{m_0}R^{m_1}L^{m_2} \cdots R^{m_{j}}$.
The sequence $m_0,....m_j$ is the same as the sequence of the
$a_i$'s in the continued fraction expansion (\cite{Ser1}).

\subsection{Connections: continued fraction expansions,  Farey words and algorithmic words}
In \cite{GK} it was established that the primitive exponents, the
F-sequence, the continued fraction expansion and the left-right
sequences are all related. The F-sequence, the left-right sequence
and the continued fraction expansion are all essentially the same.
Let $G=\langle A,B\rangle$ and $(C,D)$ be a pair of primitive words
in $A$ and $B$. Then $C$ and $D$ are {\sl primitive associates} if
$G=\langle C, D\rangle$.

\begin{theorem} [Gilman-Keen]
If the F-sequence of the algorithm is $[n_1,....,n_k]$,  let
${\frac{p}{q}}$ be the rational number with continued fraction
expansion $[n_1,...,n_k]$ and ${\frac{r}{s}}$ the rational number
with continued fraction expansion $[n_1,...,n_{k-1}]$. Let $(C,D)=
(A_k,B_k)$ be the pair of generators at which the algorithm stops.
Then $C$ is the Farey word $W_{\frac{r}{s}}(A,B)$ and $D$ is the
Farey word $W_{\frac{p}{q}}(A,B)$.
\end{theorem}

\subsection{Connections: computational complexity of the algorithm and
short geodesics}

 The computational complexity
of the $PSL(2,\mathbb{R})$ algorithm depends upon  the $F$-sequence
\cite{G21, YC}. The size of the input is measured by the trace of
the initial generators. The polynomial bound on complexity of the
algorithm was established by YC Jiang \cite{YC} and the proof
depended upon showing the fact that while Fibonacci steps made the
length of the words considered by the algorithm grow exponentially,
in that case the traces of the words considered decreased
logarithmically (as opposed to linearly).

When the group is discrete, words in the $PSL(2,\mathbb{R})$
algorithm correspond to geodesics on the surface.
 Steps in the algorithm
can be thought of as {\sl unwinding one curve about another}. A is
{\sl simple} geodesic has no self-intersections. At step $k$, the
words $(A_k, B_k)$ represent curves on the surface with fewer
self-intersections than the words $(A_{k-1},B_{k-1})$. The
F-sequence of the algorithm gives its computational complexity and
it also can be used to assign a {\sl complexity} to the curves on
the surface that the algorithm visits. The {\sl complexity} of a
curve is the number of essential self-intersections and roughly
speaking, as these decrease the length of the geodesic decreases.

The axes of $A$ and $B$ are disjoint when the
 when the trace of the multiplicative
commutator of $A$ and $B$ is positive. For this  case
\begin{theorem} [Gilman-Keen] When the algorithm stops saying that $G$ is
discrete, the algorithm has found the three shortest simple
geodesics on the surface. These are unique.
\end{theorem}

This has potential applications to the CARS program where symbolic
computation is fastest for short geodesics. The length of a geodesic
is computed through the absolute value of the trace of the
corresponding matrix, which is why the concept of trace minimizing
has an important role.

\subsection{Connection: good words, algorithm words and Farey words}

We outline another  connection between the semi-group of good words
under composition of  good words $\gs$ and addition of Farey words
$\fs$. More details about this will appear in \cite{Gconn}. First of
all we observe that for {\sl any} two words in $A$ and $B$,
$W_1(A,B)$ and $W_2(A,B)$, the operation $\gs$ is well defined:
$W(A,B) = W_2(W_1(A,B),B) = W_2 \gs W_1$. We, therefore, refer to
this operation as the {\sl good} product as opposed to the {\sl good
word} product. Any particular family of words may or may not be
closed under this product.

Assume that $W_1(X,Y)$ and $W_2(X,Y)$ are Farey words in the
generators $X$ and $Y$ with  $W_1(X,Y)= W_{\frac{p}{q}}(X,Y)$ and
$W_2(X,Y) = W_{\frac{r}{s}}(X,Y)$ where ${\frac{p}{q}}=
[n_1,...,n_j]$  and ${\frac{r}{s}}= [m_1,...,m_i]$, so that $W_1$
has $F$-sequence $[n_1,...,n_j]$ and $W_2$ has $F$-sequence
$[m_1,...,m_i]$. Then $W_1= W_{\frac{p}{q}}$ and $W_2=
W_{\frac{r}{s}}$ will not be Farey neighbors unless $|ps-qr|=1$.
However, the following is true.

\begin{prop}  Let $W(A,B)$ be the word in $A$ and $B$ with
$F$-sequence $[n_1,....,n_j,m_1,...,m_i]$. Then $W$  is the good
product $W_2 \gs W_1$.

\end{prop}

\begin{proof}
We see from corollary 2.8 and proposition 2.9 of \cite{GK} that if
the algorithm words $W_1(A,B)$ and $W_2(A,B)$, respectively, have
$F$-sequences $[n_1,...,n_j]$ and $[m_1,...,m_i]$, respectively,
then the word with $F$ sequence  $[n_1,...,n_j,m_1,...,m_i]$ is
precisely $W_2(W_1(A,B),B)$. Alternately, to see this, we begin with
$W_1$ and locate its first occurrence in the Farey diagram. We then
behave as though $W_1(A,B)$ and $B$ were the initial generators and
follow the $F$ sequence for $W_2$ down in the Farey diagram for an
initial generating pair $(W_1(A,B), B)$. The word we end at is $W$
and one can check the Farey expansion and the good word
multiplication formula.
\end{proof}

We also note that algorithm or Farey words in $ABA^{-1}$ and $B$
correspond to good words in $A$ and $B$ as follows:
\begin{prop} Let $C= ABA^{-1}$ and $D=B$.
Let $W_F$ be the algorithm/Farey word in $C$ and $D$ with primitive
exponents $(v_1,...,v_t)$ so that $W_F(C,D) = AB^{-1}A^{-1}\cdot
B^{v_1}\cdots AB^{-1}A^{-1}\cdot B^{v_t}$. Then $W_F \cdot A$ will
be a good word in $A$ and $B$.
\end{prop}
\begin{proof} Look at the exponents.
\end{proof}

\end{document}